\documentclass[12pt]{article}

\usepackage{amssymb}

\newcommand{\ud}{\mathrm{d}}
\newtheorem{theorem}{Theorem}
\newtheorem{lemma}{Lemma}
\newcommand{\bean}{\begin{eqnarray*}}
\newcommand{\eean}{\end{eqnarray*}}
\newcommand{\bea}{\begin{eqnarray}}
\newcommand{\eea}{\end{eqnarray}}
\newcommand{\eqd}{\stackrel{{\rm d}}{=}}

\newcommand{\rem}{\noindent \textbf{Remark. }}
\newcommand{\rems}{\noindent \textbf{Remarks. }}
\newcommand{\proof}{\noindent \textbf{Proof. }}
\def\bbE{{\mathbb{E}}}
\def\Exp{E}
\def\bbP{{\mathbb{P}}}
\def\bbQ{{\mathbb{Q}}}
\def\Pr{P}
\def\Var{{\mathrm{Var}}}
\def\R{\mathbb{R}}
\def\Z{\mathbb{Z}}
\def\F{{\mathcal{F}}}
\def\N{\mathbb{N}}
\def\NN{{\mathcal{N}}}
\def\eps{{\varepsilon}}

\begin{document}

\title{Logarithmic speeds for one-dimensional perturbed random walk in
random environment}

\author{M.V.~Menshikov\footnote{e-mail: \texttt{mikhail.menshikov@durham.ac.uk}}
 \\
 \normalsize
 Department of Mathematical Sciences,
 University of Durham,\\
\normalsize
 South Road, Durham DH1 3LE, England.
\and Andrew R.~Wade\footnote{e-mail: \texttt{Andrew.Wade@bris.ac.uk}}
\\
\normalsize
 Department of Mathematics,
 University of Bristol,\\
\normalsize
 University Walk, Bristol BS8 1TW, England.}

\date{May 2007}

\maketitle 

\begin{abstract}
We study the random walk in random environment on
$\Z^+=\{0,1,2,\ldots\}$, where the environment is subject to a
vanishing (random) perturbation. The two particular cases that we
consider are: (i) random walk in random environment perturbed from
Sinai's regime; (ii) simple random walk with random perturbation.
We give almost sure results on how far the random walker is
from the origin, for
almost every environment. We give both upper and
lower almost sure bounds. These bounds are of order $(\log
t)^\beta$, for $\beta \in (1,\infty)$,  depending on the
perturbation. In addition, in the ergodic cases, we give results
on the rate of decay of the stationary distribution.
\end{abstract}

\vskip 2mm


\vskip 2mm


\section{Introduction}
\label{int}

The random walk in one-dimensional random environment
in Sinai's regime (which we describe
in detail below) is a famous example of a random walk with
`logarithmic speed': after a long time $t$, the random walk is, roughly speaking,
about $(\log{t})^2$ from the origin. In this paper we give other examples of
random walks in random environments with logarithmic speeds; in these cases
the environment is subject to a random perturbation. 

Our results cover both recurrent and transient cases. In the models that
we consider,
the speed is, roughly speaking, of order $(\log{t})^\beta$, where $\beta$
depends upon the model. We shall see that for the models we consider, all $\beta \in (1,\infty)$
are attained. Examples of logarithmic transience
for random walks
(such as given in our Theorem \ref{thm8} below) are seemingly rare. The terminology `speed' is perhaps
more natural in the transient case; in the recurrent case `speed' can be thought of
as a measure of the rate of growth of the upper envelope of the random walk. Before we give our main results,
we describe the probabilistic setting in which we work.

Given an infinite sequence $\omega = (p_0,p_1,p_2,\ldots)$ such
that, for some $\delta \in (0,1/2)$, $\delta \leq p_i \leq 1-\delta$ for
all $i \in \Z^+ :=
\{0,1,2,\ldots\}$, we consider $(\eta_t(\omega); t \in
\Z^+)$ 
the nearest-neighbour random walk on $\Z^+$ defined as follows. Set $\eta_0(\omega)=r$ for
some $r \in \Z^+$, and for $n\in \N := \{1,2,\ldots\}$, \bea \label{1006bb}
 \Pr [ \eta_{t+1}(\omega) = n-1 | \eta_t(\omega) = n] & = & p_n,\nonumber\\
 \Pr [ \eta_{t+1}(\omega) = n+1 | \eta_t(\omega) = n] & = & 1-p_n =:q_n,
\eea
 and $\Pr [ \eta_{t+1}(\omega)=0 | \eta_t(\omega)=0] = p_0$, $\Pr [ \eta_{t+1}(\omega)=1 |
\eta_t(\omega)=0] = 1-p_0 =: q_0$.
The given
 form for the reflection at the origin
ensures aperiodicity, which eases some technical complications.

We call the sequence of jump probabilities $\omega$ our {\em environment}. As an example, the case $p_i=1/2$ for all $i$
gives the symmetric simple random walk on $\Z^+$.

Here, we take $\omega$ itself to be random --- then
$\eta_t(\omega)$ is a {\em random walk in random environment}
(RWRE). More precisely, $p_0,p_1,\ldots$ will be random variables
on a probability space $(\Omega, \F,\bbP)$. We describe our
particular model at the end of this section.
 The RWRE
was first studied
by Kozlov \cite{kozlov} and Solomon
\cite{solomon} (in the case where $p_i$, $i \geq 0$ form an
i.i.d.~sequence). There has been considerable interest in the RWRE recently; see for
example \cite{rev} or \cite{zeit} for surveys. 
Some authors (e.g.~\cite{sinai}) consider the RWRE with state space the whole of $\Z$. For our model
we take the case of $\Z^+$, which gives rise to a richer set of models in the sense that we can
obtain positive-recurrent behaviour.

An important case in which the random environment is homogeneous and in some sense critical
 is the so-called {\em Sinai's regime}. Here $(p_0,p_1,p_2,\ldots)$ is
a sequence of i.i.d.~random variables satisfying the condition
$\bbE[\log(p_1/q_1)]=0$, where $\bbE$ is expectation under $\bbP$.
In this case, a result dating back to Solomon \cite{solomon} says
that $\eta_t(\omega)$ is null-recurrent for $\bbP$-almost every $\omega$. Solomon's result shows that
Sinai's regime is critical
in the sense that, for an i.i.d.~random environment, $\eta_t(\omega)$
is respectively ergodic (that is positive-recurrent, here) or transient
 as $\bbE[\log(p_1/q_1)]>0$ or $\bbE[\log(p_1/q_1)]<0$.

A notable property of the RWRE in Sinai's regime is its {\em speed} --- roughly speaking
 $\eta_t(\omega)$ is of order $(\log{t})^2$ for large $t$. One way to state this more
 precisely (for another, see the discussion in Section \ref{sec3})
is in terms of `almost sure' behaviour, i.e.~results that hold $\Pr$-almost surely (a.s.) for
$\bbP$-almost every (a.e.) $\omega$. (For the remainder of this
paper, we omit the $\Pr$ and $\bbP$ when the context is clear.) This is the kind of
result we give in the present paper. In
Sinai's regime for the RWRE on $\Z^+$, almost sure upper and lower
bounds were given by Deheuvels and R\'ev\'esz (\cite{dere}, Theorem
4 in particular). A similar upper bound result was given by
Comets, Menshikov and Popov (see \cite{cmp}, Theorem 3.2), proved
 via a martingale technique related to some of the ideas in the present
paper. Sharp results are given by Hu and Shi in \cite{hushi}. In
particular, the following iterated logarithm result follows from
Theorem 1.3 of \cite{hushi}.
\begin{theorem}
\label{thmb}
 \cite{hushi} Suppose that $(p_0,p_1,p_2,\ldots)$ is
an i.i.d.~sequence with $\bbE[\log(p_1/q_1)]=0$ and $\Var(p_1) >0$.
Then
there exists a constant
 $K \in (0,\infty)$ (given
explicitly in \cite{hushi}) such that, for a.e.~$\omega$, a.s.,
 for any $\eps>0$,
\begin{itemize}
\item[(i)]
for all but finitely many $t$
\[ \frac{\eta_t(\omega)}{(\log t)^2} \leq
 (1+\eps) K  \log \log \log t; \]
\item[(ii)]
for infinitely many $t$
\[ \frac{\eta_t(\omega)}{(\log t)^2} \geq (1-\eps) K  \log \log \log t. \]
\end{itemize}
\end{theorem}
Note that `a.e.~$\omega$' is short for `$\bbP$-almost every environment $\omega$', and
`a.s.' is short for `$\Pr$-almost surely'. We use this shorthand
in the statements of all our results.
Our methods do not enable our results to be as sharp as those
in \cite{hushi}; the best that 
we obtain in Sinai's regime is included in
Theorem
\ref{thm3} below. However, we obtain a much wider array of results.

We remark that a range of {\em polynomial} speeds can be
attained in certain transient homogeneous random environment regimes (see e.g.~\cite{kks}).
In this paper we are interested primarily in {\em logarithmic}
speed results like Theorem \ref{thmb}, for random environments that are {\em asymptotically}
homogeneous.
Our main results are almost
sure upper bounds for $\eta_t(\omega)$ that are valid for a.e.~$\omega$ and all but finitely many $t$, and almost sure
lower bounds for $\eta_t(\omega)$ that are valid for a.e.~$\omega$ {\em either} for all but finitely many $t$ (if
$\eta_t(\omega)$ is transient, see e.g.~Theorem \ref{thm8}) {\em or} for infinitely many $t$
(if $\eta_t(\omega)$ is recurrent, see e.g.~Theorem \ref{thm5}).
These bounds are all of size $(\log t)^\beta$, for some $\beta \in (1,\infty)$
that is a function of $\alpha$ (the size of the perturbation),
depending on the model in question, with higher order logarithmic
corrections. 

We study two particular cases of random
environment. In the first, our environment will be a perturbation
of the i.i.d.~environment of
Sinai's regime (see Section \ref{secrwre}). In the second,
our environment will be a random perturbation of the 
simple symmetric random walk (see Section \ref{sec2}). By
studying a range of perturbations, we obtain a spectrum of
possible behaviour.

The related paper \cite{mw} employs the method of Lyapunov
functions (see \cite{fmm}) to give qualitative characteristics 
for these models (amongst somewhat more general results): specifically, 
criteria for recurrence,
transience and positive-recurrence (ergodicity, here). In the present paper
we are concerned with corresponding {\em quantitative} behaviour: specifically,
speeds (for those cases with logarithmic speed) and, secondarily,
 the rates of
decay of the stationary distribution in the ergodic cases
identified in \cite{mw}. We summarize the relevant results from \cite{mw}
at convenient points in Section \ref{results} below.

The proofs of the main results in the present paper proceeds
by relating the position
 of the random walk to some expected hitting times. The latter are analyzed (over all environments) using
 estimates for sums of independent random variables; this relies on
 (mostly well-known) strong limit theorems.

We now give a formal description of the RWRE model that we study here.
Fix $\delta \in (0,1/2)$.
Let $(\xi_i,Y_i)$, $i \in \N$, be
a sequence of i.i.d.~random vectors on some probability space $(\Omega,\F,\bbP)$, such that
\bea
\label{ue}
 \bbP [ \delta \leq \xi_1 \leq 1-\delta ] =1, \eea
 and $Y_1$ 
takes values in $[-1,1]$. The condition (\ref{ue}) is sometimes referred to as
{\em uniform ellipticity}. 
Note that we allow $Y_1$ and $\xi_1$ to be dependent. 

We fix $\alpha > 0$. For a particular realization of the
sequence $(\xi_i,Y_i)$, $i\in\N$,
we define $p_0=q_0=1/2$ and
the quantities $p_n$ and $q_n$, $n=1,2,3,\ldots$ as follows:
\bea
\label{1006b}
 p_n & := & \left\{ \begin{array}{ll}
  \xi_n +Y_n n^{-\alpha}  & ~~~~{\rm if}~~~ (\delta/2) \leq \xi_n +Y_n n^{-\alpha} \leq 1-(\delta/2) \\
  \delta/2  & ~~~~{\rm if}~~~ \xi_n +Y_n n^{-\alpha} < (\delta/2) \\
  1-(\delta/2)  & ~~~~{\rm if}~~~ \xi_n +Y_n n^{-\alpha} > 1-(\delta/2) \end{array} \right.
 \nonumber\\
 q_n & := & 1-p_n.
\eea
A particular realization of $(p_n; n\in \N)$ specifies our
random environment $\omega$,
and is given in terms of the $\xi_i$ and $Y_i$ as in (\ref{1006b}).
For a given environment $\omega$,
 the stochastic process $(\eta_t(\omega);t \in \Z^+)$ as defined at (\ref{1006bb}) is
 an irreducible, aperiodic Markov chain (under $\Pr$);
the probability measure $\Pr$ in (\ref{1006bb})
is known as the {\em quenched} measure (the measure
given a fixed environment $\omega$).

Under condition (\ref{ue}), we have that there exists $n_0 \in \N$ such that,
 for a.e.~$\omega$,
$(\delta/2) < \xi_n +Y_n n^{-\alpha} < 1-(\delta/2)$ for all $n \geq n_0$ (since the $Y_n$ are bounded). Thus,
 for all $n \geq n_0$, (\ref{1006b}) implies that, for a.e.~$\omega$,
\bean
p_n = \xi_n +Y_n n^{-\alpha}, ~~~ q_n = 1-\xi_n - Y_n n^{-\alpha},
~~~ (n \geq n_0). \eean
The  conditions on the variables
in (\ref{1006b}) ensure that, for a.e.~$\omega$, $(\delta/2) \leq p_n \leq 1 -(\delta/2)$
for all $n$ so that $p_n$ and
$q_n$ are true probabilities bounded strictly away from $0$ and
$1$, as required by our condition on $\omega$ given just before (\ref{1006bb}).

\section{Main results}
\label{results}

In this section we describe in detail two particular cases of the model formulated in
the previous section, along with our main results in each case. Then in Section
\ref{sec3} we  make
 further remarks and state some open problems.

\subsection{Perturbation of random walk in random environment in Sinai's regime}
\label{secrwre}

Now we describe our first particular case of the model given in
Section \ref{int}. For $n \in \N$ set \bea \label{0520b}
\zeta_n := \log \left( \frac{\xi_n}{1-\xi_n} \right), ~~~ Z_n :=
\frac{Y_n}{\xi_n (1-\xi_n)}. \eea With $\bbE$ denoting expectation
under $\bbP$, suppose that $\bbE[\zeta_1]=0$ and $\Var[\zeta_1]>0$
(so our environment is truly random). 
In order to formulate our results, we introduce some more notation. Set
\bea
\label{0520bx}
\lambda := \bbE [ Z_1 ] ,\eea
and also let
\bea \label{0210f} s^2:= \Var[\xi_1], ~~~ \sigma^2 := \Var[Y_1].
\eea
Under our boundedness conditions on $\xi_1$ and $Y_1$, we have
$s^2 <\infty$ and $\sigma^2 <\infty$, and under condition (\ref{ue}) we have,
$\bbP$-a.s.,
\bean
-\infty < \frac{-1}{\delta^2} \leq Z_1 \leq \frac{1}{\delta^2} < \infty.\eean

This model was introduced in \cite{mw} in somewhat more
generality, and criteria for transience, recurrence and ergodicity
given (see Theorems 6, 7 of \cite{mw}). In this case, the random
environment described in (\ref{1006b}) corresponds to a
perturbation of Sinai's regime, in the sense that, in the limit as
$n \to \infty$, we have $\bbE[\log(p_n/q_n)] \to 0$. Despite this, 
the behaviour of this model may be strikingly
different to that of Sinai's RWRE
(as 
demonstrated by our results below and also those  in
\cite{mw}), and depends on the sign of $\lambda$ as defined at
(\ref{0520bx}) (the average direction of the perturbation),
 and $\alpha$ (the size of the perturbation).

For the following results, with the definitions at (\ref{0520b}) and (\ref{0210f}),
we take $s^2>0$, $\bbE[\zeta_1]=0$,
and $\sigma^2 \geq 0$ (so, for example,
we permit the case $\bbP[Y_1=b]=1$ for some $b \in [-1,1]$, i.e.~a non-random perturbation
of Sinai's RWRE).
Of separate interest are the cases $\lambda=0$ and $\lambda \neq 0$ (where $\lambda$ is given by
(\ref{0520bx})).
The case of most interest to us here is $\lambda \neq 0$, for which the perturbation
is on average either towards $0$ ($\lambda>0$) or away from $0$ ($\lambda<0$); this
includes the case of a non-random perturbation of Sinai's RWRE. It was shown in \cite{mw} that the critical size
of the perturbation is $\alpha=1/2$: for $\alpha<1/2$ the perturbation is large enough
to disturb the null-recurrent behaviour; for $\alpha \geq 1/2$ it is too small.
By Theorem 6 of \cite{mw}, we have that if $\lambda<0$ and $\alpha<1/2$ then
$\eta_t(\omega)$ is transient for a.e.~$\omega$; if $\alpha \geq
1/2$ and $\lambda \neq 0$ then $\eta_t(\omega)$ is null-recurrent
for a.e.~$\omega$; if $\lambda>0$ and $\alpha< 1/2$ then
$\eta_t(\omega)$ is ergodic for a.e.~$\omega$.

We obtain logarithmic speeds for the $\lambda \neq 0$ case,
for the null-recurrent (Theorem \ref{thm5}), 
transient (Theorem \ref{thm8}), and ergodic (Theorem \ref{thm21}) regimes.
In the case $\lambda =0$, the critical exponent for $\alpha$ of $1/2$ is {\em
decreased}, depending on certain higher order analogues of
$\lambda$ (see the remark after Theorem 7 of \cite{mw}). Here, of the
$\lambda=0$ cases, we will only be concerned (see Theorem \ref{thm3}, below) with the special case
where $Y_1/\xi_1 \eqd -Y_1/(1-\xi_1)$, for which $\lambda=0$ and
$\eta_t(\omega)$ is null-recurrent for a.e.~$\omega$ for {\em any}
$\alpha>0$ (see \cite{mw}, Theorem 5). (Here and subsequently $\eqd$ stands for
equality in distribution.) This case is of interest because, despite the
presence of a (potentially strong) perturbation, the random walk remains null-recurrent; we
show it has logarithmic speed.

Our first result is Theorem \ref{thm5} below, which deals with
the $\lambda \neq 0$, $\alpha \geq 1/2$ case, for which $\eta_t(\omega)$
is null-recurrent for a.e.~$\omega$ (see above).
Recall
the definitions of $\lambda$, $s^2$ and $\sigma^2$ from (\ref{0520bx})
and (\ref{0210f}).

\begin{theorem}
\label{thm5} Suppose $\bbE[\zeta_1]=0$, $s^2 \in (0,\infty)$, $\lambda
\neq 0$ and $\sigma^2 \in [0,\infty)$.
\begin{itemize}
\item[(i)] Suppose $\alpha >1/2$. Then, for a.e.~$\omega$, for any $\eps>0$ we have, a.s.,
\bea 
\label{0210ac} 
0 \leq \frac{\eta_t(\omega)}{ (\log t)^2} <
(\log \log{t})^{2 +\eps},  \eea for all but finitely many $t$.
\item[(ii)]  Suppose $\alpha=1/2$. Then, for a.e.~$\omega$, for any $\eps>0$ we have, a.s.,
\bea 
\label{0210acd} 
0 \leq \frac{\eta_t(\omega)}{ (\log t )^2} <
(\log \log{t})^{4+\eps}, \eea for all but finitely many $t$.
\item[(iii)] On the other hand, for $\alpha \geq 1/2$, for a.e.~$\omega$, for any $\eps>0$ we have, a.s.,
 \bea
 \label{0920a} 
\frac{\eta_t(\omega)}{(\log{t})^2} >
(\log\log\log{t})^{-1-\eps}, \eea for infinitely many $t$.
\end{itemize}
\end{theorem}

Our next result deals with the transient case when $\lambda<0$ and $\alpha \in (0,1/2)$, and gives
a reasonably tight envelope which the random walk leaves only finitely often.
Although the random walk is transient, it is very slow: we have a striking
example of logarithmic
transience. 

\begin{theorem}
\label{thm8} Suppose $\bbE[\zeta_1]=0$, $s^2 \in (0,\infty)$, $\lambda <
0$, $\sigma^2 \in [0,\infty)$, and $\alpha \in (0,1/2)$. For
a.e.~$\omega$,
 for any $\eps>0$, we have, a.s.,
  \bea 
  \label{0427c}
 (\log  \log t)^{-(1/\alpha)-\eps} <
\frac{\eta_t(\omega)}{(\log t)^{1/\alpha}} < (\log \log
t)^{(2/\alpha)+\eps}, \eea for all but finitely many $t$. 
\end{theorem}

A case of secondary interest is that in which $Y_1/\xi_1 \eqd
-Y_1/(1-\xi_1)$. Here $\lambda=0$, and further, $\eta_t(\omega)$
is null-recurrent for a.e.~$\omega$, for any $\alpha>0$ (see
Theorem 5 of \cite{mw}). Our next result, Theorem \ref{thm3}
below, deals with this case. 

The condition $Y_1/\xi_1 \eqd
-Y_1/(1-\xi_1)$ ensures that although the perturbation may be strong, 
(roughly speaking) it balances out overall with equal strength to the left and
to the right. This intuition is supported by the fact that the random walk
remains null-recurrent. Also, included is the case
$\bbP[Y_1=0]=1$ and $\sigma^2=0$, i.e.~Sinai's regime. Thus, for our purposes,
there is no distinction between the behaviour of the RWRE perturbed from Sinai's regime
under condition
 $Y_1/\xi_1 \eqd
-Y_1/(1-\xi_1)$ and that of the RWRE in Sinai's regime itself. 
 
\begin{theorem}
\label{thm3} Suppose $\bbE[\zeta_1]=0$, $s^2 \in (0,\infty)$, $Y_1/\xi_1
\eqd -Y_1/(1-\xi_1)$,
$\sigma^2 \in [0,\infty)$,
and $\alpha >0$. \begin{itemize} \item[(i)] For a.e.~$\omega$,
 for any $\eps>0$ we have that, a.s.,
  \bean 
  0 \leq \frac{\eta_t(\omega)}{(\log
t)^2} \leq (\log \log t)^{2+\eps},  \eean
for all but finitely many $t$.
\item[(ii)] On the other hand, for a.e.~$\omega$,
for any $\eps>0$ we have that, a.s., \bean
\frac{\eta_t(\omega)}{(\log t)^2} \geq (\log \log \log
t)^{-1-\eps}, \eean
 for infinitely many $t$.
\end{itemize}
\end{theorem}

\rems (a) 
In the case of Sinai's regime ($\bbP[Y_1=0]=1$, $\sigma^2=0$), 
Theorem \ref{thm3}(i) gives similar   bounds
to \cite{dere,cmp}, but by comparison to Theorem \ref{thmb} (due to Hu and Shi \cite{hushi}),
none  of the bounds in Theorem \ref{thm3}  is particularly sharp. 
\\
(b) In the null-recurrent regimes $\lambda \neq 0$, $\alpha \geq 1/2$ (Theorem
\ref{thm5}) and $Y_1/\xi_1 \eqd
-Y_1/(1-\xi_1)$ (Theorem \ref{thm3}) we see that the position of the random walk is
essentially of order $(\log t)^2$, as in Sinai's regime (which is included
in Theorem \ref{thm3}). Thus provided we have null-recurrence we have 
the same speed. On the other hand, in the
transient case $\lambda<0$, $\alpha<1/2$ (Theorem \ref{thm8}),
the $1/\alpha$ exponent in the speed of transience is in $(2,\infty)$.
Thus for $\alpha$ increasingly small (i.e.~a
stronger perturbation), the speed increases (but is still `slow', i.e.~logarithmic).\\

In the ergodic
situations, in addition to our results on the speed of the random walk, in
the present paper we also give results on the rate of decay of the
stationary distribution $(\pi_n)$, $n\in \Z^+$,
of the Markov chain $\eta_t(\omega)$. Some analogous results for non-random environments are
given in \cite{mp}.
Theorems \ref{thm21} and \ref{thm9} below deal with the ergodic case when 
$\lambda>0$ and $\alpha \in(0,1/2)$.

\begin{theorem}
\label{thm21} Suppose $\bbE[\zeta_1]=0$, $s^2 \in (0,\infty)$, $\lambda >
0$, $\sigma^2 \in [0,\infty)$, and
$\alpha \in (0,1/2)$.  For a.e.~$\omega$, for any $\eps>0$, a.s.,
 \bean
\eta_t (\omega) \leq (1+\eps)  \left( \frac{1-\alpha}{\lambda} \right)^{1/(1-\alpha)} 
(\log t)^{1/(1-\alpha)}, 
\eean
for all but finitely many $t$, and
 \bean
\eta_t (\omega) \geq (1-\eps) \left( \frac{1-\alpha}{\lambda} \right)^{1/(1-\alpha)} (\log t)^{1/(1-\alpha)}, 
\eean
for infinitely many $t$. 
\end{theorem}

For $\alpha \in (0,1/2)$, $1/(1-\alpha) \in (1,2)$: this is `slower' than Sinai's regime.

\begin{theorem}
\label{thm9} Suppose $\bbE[\zeta_1]=0$, $s^2 \in (0,\infty)$, $\lambda >
0$, $\sigma^2 \in [0,\infty)$,
and  $\alpha \in (0,1/2)$. For a.e.~$\omega$, as $n \to \infty$
 \bea 
 \label{0922a}
\pi_n = \exp \left( -\left(\frac{\lambda}{1-\alpha}\right) n^{1-\alpha} [1+o(1)] \right).\eea
\end{theorem}

\subsection{Simple random walk with random perturbation}
\label{sec2}

Our second model again fits into the framework of (\ref{1006b}) above,
but we now take $\bbP[\xi_1=1/2]=1$ and
$\sigma^2:=\Var[Y_1]>0$. That is, we have a random perturbation of the
 symmetric simple random walk (SRW). In this case, from
(\ref{1006b}), we have $p_0=q_0=1/2$ and
for $n \in \N$
 \bea
\label{1006bc}
 p_n & := & \left\{ \begin{array}{ll}
  \frac{1}{2} +Y_n n^{-\alpha}  & ~~~~{\rm if}~~~ (\delta/2)
   \leq \frac{1}{2} +Y_n n^{-\alpha} \leq 1-(\delta/2) \\
  \delta/2  & ~~~~{\rm if}~~~ \frac{1}{2} +Y_n n^{-\alpha} < (\delta/2) \\
  1-(\delta/2)  & ~~~~{\rm if}~~~ \frac{1}{2} +Y_n n^{-\alpha} > 1-(\delta/2) \end{array} \right.
 \nonumber\\
 q_n & := & 1-p_n.
\eea
Since the $Y_n$ are bounded, we have that there exists $n_0 \in \N$ such that
 for a.e.~$\omega$ we have
$(\delta/2) < \frac{1}{2} +Y_n n^{-\alpha} < 1-(\delta/2)$ for all $n \geq n_0$. Thus,
 for a.e.~$\omega$,
 (\ref{1006bc}) implies that for all $n \geq n_0$
\bea
\label{1003a}
p_n = \frac{1}{2} +Y_n n^{-\alpha}, ~~~ q_n =
\frac{1}{2} - Y_n n^{-\alpha}, ~~~ (n \geq n_0).
\eea
The conditions
on the variables in (\ref{1006bc})
ensure that, for a.e.~$\omega$, $(\delta/2) \leq p_n \leq 1 -(\delta/2)$
 for all $n$ so that for  
$p_n$ and $q_n$ are  
bounded strictly away from $0$ and  $1$.
We see that, for a.e.~$\omega$,
 $(p_n,q_n) \to (1/2,1/2)$ as $n \to \infty$.
Thus in the limit $n \to \infty$, we coincide with
 the symmetric SRW on
$\Z^+$.

Here we do not study the case $\Var[Y_1]=\sigma^2=0$, in which  
we have a non-random perturbation of the SRW. This is an example
of the so-called Lamperti problem after \cite{lamp1} (see also \cite{harris}); for
recurrence/transience criteria 
see \cite{lamp1,mai} and Theorem 2 of \cite{mw}. From now on we assume $\Var[Y_1]=\sigma^2>0$.

The transience and recurrence properties of the model given by (\ref{1006bc})
 were
analysed in \cite{mw}. From Theorem 3(iv) of \cite{mw}, we have
that in this case if $\bbE[Y_1]<0$ and $\alpha<1$ then
$\eta_t(\omega)$ is transient for a.e.~$\omega$; if $\alpha>1$ and
$\bbE[Y_1] \neq 0$ then $\eta_t(\omega)$ is null-recurrent for
a.e.~$\omega$; if $\bbE[Y_1]>0$ and $\alpha<1$ then
$\eta_t(\omega)$ is ergodic for a.e.~$\omega$. Thus, in contrast to
the perturbation of the  {\em random} environment (as in Section \ref{secrwre}),
the critical exponent in this case is $\alpha=1$.

When $\bbE[Y_1]=0$, recurrence/transience properties depend on the
higher moments of $Y_1$ (see the remark after Theorem 3 of
\cite{mw}).
Of interest to us in the present paper is the case in which the
distribution of $Y_1$ is symmetric, that is $Y_1 \eqd -Y_1$ (and
$\bbE[Y_1]=0$). In this case (see Theorem 3(iii) of \cite{mw})
$\eta_t(\omega)$ is null-recurrent for a.e.~$\omega$, for {\em any}
$\alpha>0$. In this case we obtain our logarithmic behaviour (see
Theorem \ref{thm0}), in the domain $\alpha \in (0,1/2)$. We also obtain
logarithmic bounds in the ergodic case mentioned above (see Theorem \ref{thm20}).
 
\begin{theorem}
\label{thm0}
Suppose $\bbP[\xi_1=1/2]=1$, $Y_1 \eqd -Y_1$,
$\sigma^2 \in (0,\infty)$, $\alpha\in (0,1/2)$.
\begin{itemize}
\item[(i)] For a.e.~$\omega$, for any
$\eps>0$, a.s.,
 \bea 
 \label{0210a} 
 0 \leq
\frac{\eta_t(\omega)}{ (\log t)^{2/(1-2\alpha)}} \leq (\log \log
t)^{(2/(1-2\alpha))+\eps}, 
\eea for all but finitely many $t$. \item[(ii)] On the other hand,
for a.e.~$\omega$, for any $\eps>0$, a.s.,
 \bea 
 \label{0920c}
\frac{\eta_t(\omega)}{(\log{t})^{2/(1-2\alpha)}} \geq ( \log \log
\log t )^{-(1/(1-2\alpha))-\eps} , \eea for infinitely many $t$.
\end{itemize}
\end{theorem}
\rem Note that for $\alpha \in (0,1/2)$, $2/(1-2\alpha)$ is in $(2,\infty)$.
In the limit $\alpha \downarrow 0$, we approach Sinai's regime in the sense that, for
fixed $\omega$ and each $n$,
$(p_n,q_n) \to (\frac{1}{2}+Y_n,\frac{1}{2}-Y_n)$ where \[ \bbE \left [ \log \left(\frac{(1/2)+Y_n}{(1/2)-Y_n}
\right) \right] =
\bbE [ \log ((1/2)+Y_n) ] - \bbE [ \log ((1/2)-Y_n)]= 
 0 \] when $Y_1 \eqd -Y_1$. Thus it is not surprising that in the limit
$\alpha \downarrow 0$, Theorem \ref{thm0} approaches Theorem \ref{thm3} (which includes
Sinai's regime).\\

Theorems \ref{thm20} and \ref{thm10} below deal with the ergodic case when $\bbE[Y_1]>0$ and 
$\alpha \in(0,1)$. 
Note that when $\alpha \in (0,1)$, $1/(1-\alpha) \in (1,\infty)$.

\begin{theorem}
\label{thm20} Suppose $\bbP[\xi_1=1/2]=1$, $\bbE[Y_1] >
0$, $\sigma^2 \in (0,\infty)$, and $\alpha \in (0,1)$.  For a.e.~$\omega$, for any $\eps>0$, 
a.s.,
\bean
\eta_t (\omega) \leq (1+\eps) 
 \left( \frac{1-\alpha}{4 \bbE[Y_1]} \right)^{1/(1-\alpha)} 
(\log t)^{1/(1-\alpha)},
\eean
for all but finitely many $t$, and
 \bean
\eta_t (\omega) \geq (1-\eps) 
 \left( \frac{1-\alpha}{4 \bbE[Y_1]} \right)^{1/(1-\alpha)}
 (\log t)^{1/(1-\alpha)}, 
\eean
for infinitely many $t$. 
\end{theorem}

The next result 
gives the rate of decay of the stationary distribution $(\pi_n)$:
as in Theorem \ref{thm9}, the decay is
sub-exponential.

\begin{theorem}
\label{thm10} Suppose $\bbP[\xi_1=1/2]=1$, $\bbE[Y_1] >
0$, $\sigma^2 \in (0,\infty)$, and $\alpha \in (0,1)$.  For a.e.~$\omega$, as $n \to \infty$
 \bea \label{0922c}
\pi_n = \exp \left( -\left(\frac{4 \bbE[Y_1]}{1-\alpha}\right) n^{1-\alpha} [1+o(1)] \right).\eea
 \end{theorem}

\subsection{Further remarks and open problems}
\label{sec3}

Our results give an indication of the `almost sure' behaviour of
$\eta_t(\omega)$, and there is scope for tightening our  bounds. 
Also of interest is 
the so-called {\em annealed} behaviour of the RWRE (averaged
 over all environments). Sinai's result \cite{sinai}
for the random walk in i.i.d.~random environment on $\Z$ with $\bbE[\log(p_1/q_1)]=0$ showed
(roughly speaking) that
$\eta_t(\omega)$ divided by $(\log{t})^2$ converges in distribution to some random variable
as $t \to \infty$.
The result is stated in terms of the annealed probability measure $\bbQ$ given by
\[ \bbQ [\cdot] = \int_\Omega \Pr [\cdot] \ud \bbP [\omega].\]
Golosov \cite{gol1} showed that for the RWRE on $\Z^+$ in Sinai's regime
\[ \bbQ \left[\frac{\eta_t (\omega)}{(\log t)^2} \leq u \right] \longrightarrow F(u), ~~~ u \in \R,\]
as $t \to \infty$, where $F$ is a known distribution function. 
See also \cite{gol2,gol3,kes,cp} for
related results. The annealed behaviour of our models is also of interest. In
particular, under the conditions of Theorem \ref{thm8} do we have (analogously to the results of Sinai-Golosov
\cite{sinai,gol1}) that as $t \to \infty$
\[ \bbQ \left[ \frac{\eta_t(\omega)}{(\log t)^{1/\alpha}} \leq u \right] \longrightarrow G(u), ~~~ u \in \R,\]
for some $G$? We do not address this
question in the present paper.

One can obtain $L^p$ analogues of our results, with the methods used here (compare Theorem
3.2 of \cite{cmp}). For example, under the conditions of Theorem
\ref{thm5}, analogously to (\ref{0210ac}), 
for any $p \geq 1$, for any $\eps>0$,
for a.e.~$\omega$, as $t\to \infty$
\bean \frac{\eta_t(\omega)}{ (\log t
)^{2+\eps} }  \to 0, ~\textrm{in} ~L^p. \eean  

The methods of the present paper are well suited to logarithmic speeds, since
they are based on an analysis of the expected hitting times of the
random walk; some standard estimates using the submartingale property,
Markov's
inequality and the (first)
Borel-Cantelli lemma lead to some rather sharp
results, since these expected times are exponentially large. Of interest would be
results for the cases of the SRW
 with random perturbation
that are
not covered by the theorems of Section \ref{sec2}. For example,
if $Y_1 \eqd -Y_1$ but $\alpha > 1/2$, we expect SRW-like
behaviour. On the other hand, if $\bbE[Y_1] \neq 0$, we suspect
that  $\eta_t(\omega)$ will behave in a similar way to the Lamperti problem
mentioned above: roughly speaking, we expect
  SRW-like behaviour for $\alpha >1$, while in the transient regime
($\alpha<1$ and $\bbE[Y_1]<0$) we have $\eta_t (\omega) \sim t^{1/(1+\alpha)}$. 
Another open problem is 
the behaviour of this model when $\alpha=1$ (this case was
not covered in \cite{mw}).
We hope to address some of these
issues in future work.

\section{Preliminaries}
\label{prelim}

Before we prove our main results in Section \ref{secprfs}, we give some preparatory results.
First, in Section \ref{strng},
we present some technical lemmas concerning the behaviour of sums of independent random variables; some are well-known results, others we prove.
Then, in Section \ref{hitting},
 we give the main apparatus of our proofs, based on some hitting time results.

\subsection{Some strong theorems for sums of independent random variables}
\label{strng}

The following result is due to Sakhanenko \cite{sak1,sak2,sak3},
and is contained in Theorem A* of the more readily obtainable
paper by Shao \cite{shao}.
\begin{lemma}
\label{lem1025a}
Let $X_1,X_2, \ldots$ be independent random variables
with $\Exp[X_i]=0$, $\Var [X_i] = \sigma_i^2 \in (0,\infty)$ for
$i \in \N$. Suppose that the
$X_i$ are uniformly bounded, i.e., for some $B \in (0,\infty)$, $\Pr[|X_i| > B] =0$ for all $i$.
For $n \in \N$, set
\[ s^2_n := \sum_{i=1}^n \sigma_i^2 .\]
Then, there exists (possibly on an enlarged probability space)
a sequence of independent normal random variables $(W_1,\ldots,W_n)$ with
$\Exp[W_i]=0$, $\Var[W_i]=\sigma_i^2$ for $1 \leq i \leq n$ such that
a.s.,
\bean
\left| \sum_{i=1}^n X_i - \sum_{i=1}^n W_i \right|
 \leq \frac{1}{A} \log (s_n^2) ,
\eean
for all but finitely many $n$,  
where $A\in (0,\infty)$ is a constant.
\end{lemma}

We will need a form of the Law of the Iterated
Logarithm. The following result is a consequence of
Theorem 7 of \cite{feller}.
\begin{lemma}
\label{itlog} Let $X_1, X_2, \ldots$ be 
independent, uniformly bounded
 random variables with $\Exp[X_i]=0$,
$\Var[X_i^2]=\sigma_i^2 \in (0,\infty)$ for $i \in \N$.  For $n\in \N$, set
$s_n^2 := \sum_{i=1}^n
\sigma_i^2$.
 Suppose that $s_n \to \infty$ as $n \to \infty$.
  Then, for any $\eps>0$, a.s.,
\bean
 \left| \sum_{i=1}^n X_i \right| \leq s_n ((2+\eps)\log\log (s_n^2))^{1/2} , \eean
  for all but finitely many $n$. 
\end{lemma}

We will also need the following `inverse iterated logarithm law'
due to Hirsch (Theorem 2 of \cite{hirsch}; see also Theorem 3.1 of \cite{csaki}).

\begin{lemma}
\label{hlem}
Let $X_1,X_2,\ldots$ be i.i.d., uniformly bounded
random variables with
$\Exp[X_1]=0$, $\Var[X_1] \in (0,\infty)$. For $x \geq 0$, let $a(x)>0$ be a
nonincreasing function
such that $x^{1/2}a(x)$ is eventually increasing and
\bea
\label{sumn}
\sum_{n=1}^\infty \frac{a(n)}{n}  < \infty.\eea
Then, a.s., 
\bean
\max_{1 \leq i \leq n} \sum_{j=1}^i X_j \geq n^{1/2} a(n),
\eean
for all but finitely many $n$.
\end{lemma}

We will also need the following extension of part of Hirsch's result to independent
non-identically
distributed random variables.
\begin{lemma}
\label{lem0601} Let $X_1,X_2,\ldots$ be 
independent,
 uniformly bounded
random variables with $\Exp[X_i]=0$,
$\Var[X_i]=\sigma_i^2$ for $i\in\N$, where
$0<\sigma_i^2<M<\infty$ for all $i$. 
 Set $s_n^2:= \sum_{i=1}^n \sigma_i^2$ for
$n \in \N$. Suppose that $s_n \to \infty$ as $n \to \infty$.
For $x \geq 0$, let $a(x)>0$ be a nonincreasing function
such that $x^{1/2}a(x)$ is eventually increasing, (\ref{sumn}) holds,
and
\bea
\label{0531a}
\lim_{n \to \infty} \frac{\log{n}}{n^{1/2} a(n)} = 0 .\eea
Then, for some constant $C \in (0,\infty)$, a.s.,
\bea
\label{0210bb}
\max_{1 \leq i \leq n} \sum_{j=1}^i X_j \geq C s_n a(s_n^2),
\eea
for all but finitely many $n$.
\end{lemma}
\proof
By Lemma \ref{lem1025a}, we can redefine the $X_i$, $i\in\N$
on a richer probability space along with a sequence of independent
normal random variables $W_i$, $i\in\N$
with $\Exp[W_i]=0$ and $\Var[W_i]=\sigma_i^2$, such that, a.s.,
\[ \left| \sum_{j=1}^i X_j - \sum_{j=1}^i W_j \right|
 \leq A^{-1} \log (s_i^2) \leq C \log i ,
\]
for all but finitely many $i$, for some $A, C\in(0,\infty)$. Thus,
a.s.,
\bea
\label{0531b}
 \left| \max_{1 \leq i \leq n} \sum_{j=1}^i X_j - \max_{1 \leq i \leq n}
\sum_{j=1}^i W_j \right|
 \leq \max_{1 \leq i \leq n} \left|  \sum_{j=1}^i X_j - \sum_{j=1}^i W_j \right|
\leq C \log n,
\eea
for all but finitely many $n$.
For $n \in \N$, set
\bea
\label{0601f}
 h(n) := \min \{ m \in \N: s_m^2 \geq n \}.\eea
There exists a standard Brownian motion $(B(n);n \geq 0)$ and a
sequence of independent normal random variables $\delta_n \sim
\NN(0,s_{h(n)}^2-n)$, $n\in\N$, independent of $(B(n);n \geq 0)$, such
that
\[ B(n) + \delta_n = \sum_{i=1}^{h(n)} W_i ,\]
for each $n \in \N$.
Now,
\[ \max_{1 \leq i \leq h(n)} \sum_{j=1}^i W_j \geq
\max_{1 \leq i \leq n} \sum_{j=1}^{h(i)} W_j = \max_{1 \leq i \leq
n} \left( B(i) + \delta_i \right).\] Hence \bea \label{0601a}
 \max_{1 \leq i \leq h(n)} \sum_{j=1}^i W_j \geq \max_{1 \leq i \leq n} B(i)
- \max_{1 \leq i \leq n} \delta_i.\eea 
Since
$\Var(\delta_i)=s_{h(i)}^2-i \leq \sigma_{h(i)}^2 <M< \infty$, and
$\delta_i$, $i \in \{1,\ldots,n\}$ are independent normal random
variables, we have that, a.s., \bea
\label{0601b}
 \max_{1 \leq i \leq n} \delta_i \leq \log n ,\eea
for all but finitely many $n$ 
(this follows
from standard tail bounds on the normal distribution (see e.g.~\cite{durrett}, p.~9)
and the Borel-Cantelli lemma). Suppose that $a(\cdot)$ satisfies
the conditions of this lemma. Now, for a sequence $Y_1, Y_2,\ldots$ of i.i.d.~normal
random variables with $\Exp[Y_1]=0$ and $\Var[Y_1]=1$, we have  by Lemma \ref{hlem} that a.s., 
\bea
\label{0601c}
 \max_{1 \leq i \leq n} B(i) = \max_{1 \leq i \leq n} \sum_{j=1}^i Y_j
\geq n^{1/2} a(n) ,\eea for all but finitely many $n$. So from
(\ref{0601a}), (\ref{0601b}), (\ref{0601c}) and condition
(\ref{0531a}), a.s., \bea \label{0601d} \max_{1
\leq i \leq h(n)} \sum_{j=1}^i W_j \geq n^{1/2} a(n) - \log n
\geq C n^{1/2} a(n) ,\eea for all but finitely many $n$ and some $C \in (0,\infty)$.
Since $\sigma_i^2>0$ for all $i$,
 we have from (\ref{0601f}) that $h(s_n^2)=n$; thus by (\ref{0531b}) and (\ref{0601d})
we have that, a.s.,
 \bean \max_{1 \leq i
\leq n} \sum_{j=1}^i X_j  \geq  \max_{1 \leq i \leq n}
\sum_{j=1}^i W_j - C \log n  \geq  C' (s_n^2)^{1/2} a(s_n^2) - C \log
n ,\eean for some $C,C' \in(0,\infty)$, for all but finitely many
$n$. Then by the conditions on $s_n^2$ and $a(\cdot)$,
(\ref{0210bb}) follows. $\square$\\

The next two lemmas will be needed for some more delicate estimates (e.g.~in the proof
of Theorem \ref{thm8}) where we need to deal with certain moving sums. 
The following lemma is a corollary to a result of Hirsch \cite{hirsch}.

\begin{lemma}
\label{lbd1}
Let $X_1,X_2,\ldots$ be independent, uniformly bounded random variables with
 $\Exp[X_i]=0$, $\Var[X_i] =\sigma^2 \in (0,\infty)$ for  $i \in \N$. 
For $x \geq 1$, let $b(x)$ be a nondecreasing, integer-valued 
function such that 
 for some $\beta >0$ and $x_0 \in (0,\infty)$,
 $x^\beta \leq b(x) \leq x$ for all $x \geq x_0$.
Then
for any $\eps>0$, a.s.,
\[ \max_{1 \leq i \leq n} \max_{1 \leq j \leq b(i)}
\sum_{k=i-j+1}^i X_k \geq  (b(n/2))^{1/2} (\log n)^{-1-\eps},\]
for all but finitely many $n$.
\end{lemma}
\proof
For fixed $i$, note that
\[ \max_{1 \leq j \leq b(i)}
\sum_{k=i-j+1}^i X_k \eqd \max_{1 \leq j \leq b(i)}
\sum_{k=1}^j Y_k,\]
where $Y_1,Y_2,\ldots$ are independent random variables
with $Y_k \eqd X_{i+1-k}$ for each $k$.  Fix $\eps>0$.
Let $E_i$ denote the event 
\[ E_i := \left\{
\max_{1 \leq j \leq b(i)}
\sum_{k=i-j+1}^i X_k \leq (b(i))^{1/2} (\log b(i))^{-1-\eps}
\right\}.\]
Then
Corollary 1 of Hirsch \cite{hirsch} implies that there are
 absolute constants $C, C' \in (0,\infty)$ such that
for all $i \geq x_0$,
\[ \Pr [ E_i ]
\leq C (\log b(i))^{-1-\eps} 
\leq C' (\log i)^{-1-\eps},\]
since $b(i) \geq i^\beta$.
Consider the subsequence $i=2^m$ for $m=1,2,\ldots$. Then 
\[ \sum_{m=1}^\infty \Pr [E_{2^m}] \leq C \sum_{m=1}^\infty m^{-1-\eps} < \infty.\]
Hence by the (first) Borel-Cantelli lemma,
 a.s., there
is a finite $m_0$ (with $2^{m_0} \geq x_0$)
such that, for all $i=2^m$ with $m \geq m_0$,
\[ \max_{1 \leq j \leq b(i)}
\sum_{k=i-j+1}^i X_k \geq (b(i))^{1/2} (\log b(i))^{-1-\eps} \geq  (b(i))^{1/2} (\log i)^{-1-\eps},\]
since $b(i) \leq i$.
Each $n \geq 2$ satisfies
 $n \in [2^m,2^{m+1})$ for some $m \in \N$; then, a.s., 
\bean \max_{1 \leq i \leq n}
\max_{1 \leq j \leq b(i)}
\sum_{k=i-j+1}^i X_k 
\geq
\max_{1 \leq i \leq 2^m}
\max_{1 \leq j \leq b(i)}
\sum_{k=i-j+1}^i X_k \\
\geq
\max_{1 \leq j \leq b(2^m)}
\sum_{k=2^m-j+1}^{2^m} X_k 
\geq (b(2^m))^{1/2} (\log (2^m))^{-1-\eps},\eean
provided $m \geq m_0$.
Hence, since $n \geq 2^m > n/2$, a.s.,
\[ \max_{1 \leq i \leq n}
\max_{1 \leq j \leq b(i)}
\sum_{k=i-j+1}^i X_k \geq (b(n/2))^{1/2} (\log n)^{-1-\eps},\]
for all $n \geq 2^{m_0}$. 
$\square$

\begin{lemma}
\label{ubd1}
Let $X_1,X_2,\ldots$ be independent,
uniformly bounded
random variables with $\Exp [X_i]=0$ for all $i\in\N$.
Then
there exists $C \in (0,\infty)$ such that, a.s.,
for all but finitely many $i$
\[ \left|
\sum_{k=i-j+1}^i X_k
\right|
 \leq  C j^{1/2} (\log i)^{1/2},\]
for all $j=1,2,\ldots,i$.
\end{lemma}
\proof For fixed $i$,
$Y^i_j := \sum_{k=i-j+1}^i X_k$
is a martingale over $j=1,2,\ldots,i$, with uniformly bounded
increments. Hence the Azuma-Hoeffding inequality (see e.g.~\cite{hoef})
implies that for some $B \in (0,\infty)$, for all $j=1,\ldots,i$, for $t>0$,
\[ \Pr [ |Y^i_j| \geq t ] \leq 2 \exp ( -B^{-1} j^{-1} t^2 ).\]
Thus for a suitable   $C<\infty$, for $j \leq i$, 
$\Pr [ |Y^i_j | \geq C j^{1/2} (\log i)^{1/2} ] \leq i^{-3}$. Then 
\[ \sum_{i=1}^\infty \sum_{j=1}^i \Pr [ | Y^i_j| \geq C j^{1/2} (\log i)^{1/2} ]
\leq \sum_{i=1}^\infty i^{-2}  < \infty.\]
Hence the (first) Borel-Cantelli lemma implies that,
 a.s., there are only finitely
many pairs $(i,j)$ (with $j \leq i$) for which $|Y^i_j| \geq  C j^{1/2} (\log i)^{1/2}$.
$\square$

\subsection{Hitting times results}
\label{hitting}

For the proofs of our main results, we will use the expected
hitting times for the random walk $\eta_t (\omega)$  as defined at (\ref{1006bb}).
 For the remainder of this section,
we work in the quenched
setting (i.e.~with fixed environment $\omega = (p_0,p_1,\ldots)$ throughout).
For $0 \leq m < n$,
let $\tau_{m,n}$ denote the time when $\eta_t(\omega)$ first hits $n$, starting from $m$.
That is, with the convention $\min \emptyset = +\infty$,
 \bea
\label{1001a}
\tau_{m,n} := \min \{ t \geq 0: \eta_t(\omega) = n | \eta_0(\omega)=m \}.\eea
For our proofs in Section \ref{secprfs}, we take $\eta_0(\omega)=r=0$
for ease of exposition; the proofs easily extend to general $r \in \mathbb{Z}^+$.
For fixed $\omega$, let $T(0):=0$, and for $n\in\N$ let
$T(n):=\Exp[ \tau_{0,n}]$.
For $i=0,1,2,\ldots$,
write $\Delta_i:=T(i+1)-T(i)=\Exp[\tau_{i,i+1}]$, so that $\Delta_i$ is the expected
time taken for $\eta_t(\omega)$ to hit $i+1$, starting at $i$. Then standard arguments
yield $T(n) = \sum_{i=0}^{n-1} \Delta_i$ with
$\Delta_0=1/q_0$ and for $i \geq 1$
\[ \Delta_i = 1+ p_i ( \Delta_{i-1} + \Delta_i).\]
We then obtain the following classical result.
\begin{lemma}
\label{lemexp} Let $\omega$ be fixed. For $n \in \N$, we have that
$T(n) = \sum_{i=0}^{n-1} \Delta_i$, and for $i \geq 0$,
 $\Delta_i$
is given (with the convention that an empty product is 1) by
\bea
\label{0420b}
 \Delta_i & = & \sum_{j=0}^i q_{i-j}^{-1} \prod_{k=i-j+1}^{i}
 \frac{p_k}{q_k}
 =  \frac{1}{q_i} + \frac{p_i}{q_i q_{i-1}} + \cdots + \frac{p_i
p_{i-1} \cdots p_1}{ q_i q_{i-1} \cdots q_1 q_0} .
 \eea
\end{lemma}

The following fact will be very useful. That is,  for a fixed
environment,
 $T(\eta_t(\omega))$ is a submartingale with respect to the
 natural filtration (for a 
 closely related supermartingale, see \cite{cmp}, equation (6)).
  In particular, we have the following.
\begin{lemma}
\label{1002a}
For fixed $\omega$, any $t \in \Z^+$ and any $n\in
\Z^+$,
\bea \label{1002b} \Exp [ T( \eta_{t+1}(\omega)) -
T(\eta_t(\omega)) | \eta_t(\omega) = n ] = 1.\eea
\end{lemma}
\proof For $n \geq 1$, we have \bean & &
\Exp [ T( \eta_{t+1}(\omega))
- T(\eta_t(\omega)) | \eta_t(\omega) = n ] \\
& = & p_n (T(n-1)-T(n))
+q_n(T(n+1)-T(n)) \\
& = & q_n \Delta_n - p_n \Delta_{n-1} = 1,\eean by (\ref{0420b}).
Also,
\[ \Exp [ T( \eta_{t+1}(\omega))
- T(\eta_t(\omega)) | \eta_t(\omega) = 0 ] = q_0 T(1) = 1,\] since
$T(1)=\Delta_0=1/q_0$. $\square$\\

We can now state the result that will be our main tool in proving
almost sure upper and lower bounds for $\eta_t(\omega)$, using the
expected hitting times $T(n)$.
\begin{lemma}
\label{1002e} For a given environment $\omega$, suppose that there
exist two nonnegative, increasing, continuous functions $g$ and
$h$ such that,
\[ g(n) \leq T(n) \leq h(n) ,\]
for all $n \in \Z^+$. 
Then:
\begin{itemize}
\item[(i)] For any $\eps>0$, a.s., for all but finitely many $t$, 
 \bea \label{001}
 \eta_t(\omega) \leq  g^{-1}
( (2t)^{1+\eps} ). \eea 
\item[(ii)] 
A.s., for
infinitely many $t$,
\bea \label{003}
(\eta_t(\omega))^2 h( \eta_t(\omega)) \geq t. \eea
\end{itemize}
\end{lemma}
\rem
 In the transient case  we want to do better
 (for Theorem \ref{thm8}) than part (ii) here, to
  give a lower bound for $\eta_t(\omega)$ that holds all but finitely often.
 See the proof of Theorem \ref{thm8} below.\\

\noindent
{\bf Proof of Lemma \ref{1002e}.}
Throughout we work in fixed
environment $\omega$. First we prove part (i).
From (\ref{1002b}), we have that for any $t \in \Z^+$
\[ \Exp [ T(\eta_{t+1}(\omega)) - T(\eta_{t}(\omega))]
= \sum_{n=0}^\infty \Pr[ \eta_t(\omega)=n] = 1.\] Then, given that
$\eta_0(\omega)=0$,
for all $t \in \Z^+$ we have \bea \label{0210j}
\Exp[T(\eta_t(\omega))] = t. \eea
 To prove (\ref{001}),
we modify the idea of the proof of Theorem 3.2 of \cite{cmp}. 
Since $T(\eta_t(\omega))$ is a nonnegative submartingale (see Lemma \ref{1002a}),
 Doob's
submartingale inequality (see e.g.~\cite{williams}, p.~137)
implies that, for $t>0$, for any $\eps>0$,
 \bea \label{1002d}
 \Pr \left[ \max_{0 \leq s \leq t} T(\eta_s(\omega)) \geq
  t^{1+\eps} \right] \leq t^{-1-\eps} \Exp[ T(\eta_t(\omega)) ] =  t^{-\eps},
\eea 
using (\ref{0210j}).
Also,
 given that $T(n)
\geq g(n)$ for all $n$, we have, for $t>0$,
 \bea
 \label{1002c} \Pr \left[
\max_{0 \leq s \leq t} T(\eta_s(\omega))
\geq t^{1+\eps} \right]
 \geq  \Pr \left[ \max_{0 \leq s \leq t} g(\eta_s(\omega)) \geq t^{1+\eps} \right] \nonumber\\
  = \Pr \left[ g \left( \max_{0 \leq s \leq t} \eta_s(\omega) \right) \geq t^{1+\eps} \right] ,\eea
  since $g$ is increasing. Hence from (\ref{1002d}) and (\ref{1002c}), for $t>0$,
 \[ 
  \Pr \left[   \max_{0 \leq s \leq t} \eta_s(\omega)   \geq g^{-1} ( t^{1+\eps})  \right] 
  \leq t^{-\eps}  .\]
 Thus along the subsequence $t=2^m$ for $m=0,1,2,\ldots$, 
  the (first) Borel-Cantelli lemma implies that, a.s., the event in the last display
  occurs only finitely often, and in particular there exists $m_0 <\infty$
  such that for all $m \geq m_0$
  \[ \max_{0 \leq s \leq 2^m} \eta_s(\omega)   \leq g^{-1} ( (2^m)^{1+\eps}).\]
  Every $t$ sufficiently large has $2^m \leq t < 2^{m+1}$ for some $m \geq m_0$; then,
  a.s.,
  \[ \eta_t (\omega) \leq
  \max_{0 \leq s \leq t} \eta_s(\omega) 
  \leq  \max_{0 \leq s \leq 2^{m+1}} \eta_s(\omega) 
   \leq g^{-1} ( (2^{m+1})^{1+\eps}) ,\]
   for all but finitely many $t$. Now since $2^{m+1} \leq 2t$ and
   $g^{-1}$ is increasing, 
(\ref{001}) follows.

Now we prove part (ii).
Recall the definition of
$\tau_{0,n}$ at (\ref{1001a}). By Markov's inequality, we have
that for $n \in \N$
\[ \Pr [ \tau_{0,n} > n^2 T(n) ] = \Pr [ \tau_{0,n} > n^2 \Exp[ \tau_{0,n}] ]  \leq n^{-2} .\]
Then, by the (first) Borel-Cantelli lemma, a.s., 
$\tau_{0,n} > n^2 T(n)$ for only finitely
many $n$. Thus, given that $T(n) \leq h(n)$ for all $n$,
 we have that a.s., for all but finitely
many $n$, 
$\tau_{0,n} \leq n^2 h(n)$.

Given $\omega$, $\eta_t(\omega)$ is an irreducible
Markov chain on $\Z^+$, hence $\limsup_{t \to \infty} \eta_t (\omega) = +\infty$ a.s..
Thus a.s.~there exists an infinite subsequence of $\N$,
$t_1,t_2,t_3,\ldots$ (one can take, 
for each $i$, $t_i=\tau_{0,i}$, the time of the first visit of $\eta_t$ to $i$), 
such that $\eta_{t_i}(\omega) \to \infty$ as
$i \to \infty$. That is, a.s.,
\[ t_i \leq \eta_{t_i}(\omega)^2 h(\eta_{t_i}(\omega)) .\]
There are infinitely many such $t_i$, 
and so we have (\ref{003}).
$\square$

\section{Proofs of main results}
\label{secprfs}

To prove our main results, we employ the machinery given in the previous
section:
 we obtain,
 via the results in Section \ref{strng},
suitable
functions $g$, $h$ such that $g(n) \leq T(n) \leq h(n)$ (for a.e.~$\omega$), and then apply
 Lemma
\ref{1002e}. 

We consider $T(n)$ as given in Lemma \ref{lemexp}. Recalling 
the definition of $\Delta_i$ at (\ref{0420b}), we can write
(interpreting an empty sum as zero) for $i \geq 0$
\bea
\label{1001c}
\Delta_i  =  \sum_{j=0}^i q_{i-j}^{-1} \exp \sum_{k=i-j+1}^i
\log ( p_k / q_k )
.\eea
The following result gives general bounds on $T(n)$.

\begin{lemma}
\label{lowerbd}
For a fixed environment $\omega$, 
for all $n \geq 1$
\bea
\label{ffff}
T(n) \geq \exp \max_{1 \leq i \leq n-1} \sum_{k=1}^i \log (p_k/q_k)
, \eea
and for some $C\in(0,\infty)$, for all $n \geq 1$,
\bea
\label{eeee}
T(n) \leq
C n^2 \exp \left( \max_{0 \leq i \leq n-1} \sum_{k=1}^i \log (p_k/q_k) 
+ \max_{0 \leq i \leq n-1} \sum_{k=1}^i (-\log (p_k/q_k)) \right)
.\eea
\end{lemma}
\proof
Since a sum of nonnegative
terms is bounded below by its largest term,
\bea
\label{ab1}
 T(n) = \sum_{i=0}^{n-1} \Delta_i 
 \geq \max_{1 \leq i \leq n-1} \Delta_i
\geq \max_{1 \leq i \leq n-1} \max_{1 \leq j \leq i} \exp \sum_{k=i-j+1}^i \log (p_k/q_k) ,\eea
using (\ref{1001c}) and the fact that $q_{i-j}^{-1} \geq 1$. Now for $i \in \N$
\bea
\label{ab2}
 \max_{1 \leq j \leq i} \sum_{k=i-j+1}^i \log (p_k/q_k)
\geq \sum_{k=1}^i \log (p_k/q_k),\eea
so that by (\ref{ab1}) and (\ref{ab2}),
\[ T(n) \geq \max_{1 \leq i \leq n-1}  \exp \max_{1 \leq j \leq i} \sum_{k=i-j+1}^i \log (p_k/q_k)
\geq \max_{1 \leq i \leq n-1}  \exp  \sum_{k=1}^i \log (p_k/q_k),\]
and the lower bound in the lemma follows.

For the upper bound, we have from (\ref{1001c}) that
\bea
\label{ab10}
 T(n) \leq n \max_{0 \leq i \leq n-1} \Delta_i
\leq \delta^{-1} n(n+1) \max_{0 \leq i \leq n-1} \max_{0 \leq j \leq i} \exp \sum_{k=i-j+1}^i \log (p_k/q_k) ,\eea
since $q^{-1}_{i-j} \leq \delta^{-1}$ with $\delta$ as at (\ref{ue}). Now
\bea
\label{ab20}
 \max_{0 \leq j \leq i} \sum_{k=i-j+1}^i \log (p_k/q_k)
=  \sum_{k=1}^i \log (p_k/q_k)
+ \max_{0 \leq j \leq i} \sum_{k=1}^{i-j} (-\log (p_k/q_k)) \nonumber\\
=  \sum_{k=1}^i \log (p_k/q_k)
+ \max_{0 \leq j \leq i} \sum_{k=1}^{j} (-\log (p_k/q_k)).\eea
Thus  from (\ref{ab10}) and (\ref{ab20}), for $C \in (0,\infty)$ and all $n \geq 1$
\bean
 T(n) \leq C n^2 \exp \left( \max_{0 \leq i \leq n-1} \sum_{k=1}^i \log (p_k/q_k)
+ \max_{0 \leq i \leq n-1} \max_{0 \leq j \leq i} \sum_{k=1}^{j} (-\log (p_k/q_k)) \right) 
.\eean
Then the upper bound in the lemma follows.
 $\square$\\

We start with the proof of Theorem \ref{thm0} for expository purposes.
 The proof of Theorem \ref{thm0} will then serve
as a prototype for subsequent proofs. As previously mentioned,
we take $\eta_0(\omega)=0$ for the purposes of the proofs that follow (without loss of generality).

\subsection{Proof of Theorem \ref{thm0}}
\label{secprf1}

For fixed $\omega$, by Lemma \ref{lemexp},
the expected hitting time $T(n)$ is
 expressed
in terms of $\log (p_n/q_n)$. To prepare for the proof, we
note that under the conditions of Theorem \ref{thm0}
$p_n$ and $q_n$ have the same distribution, so
\bea
\label{0608a}
 \bbE [ \log (p_n/q_n) ] = \bbE [ \log p_n]-\bbE [\log q_n] = 0.\eea
By (\ref{1003a}), Taylor's theorem
and the boundedness of the $Y_n$,
for a.e.~$\omega$,
\bean
 \log p_n & = & 
 \log (1/2) + \log ( 1+ 2Y_n n^{-\alpha} ) \\ 
 & = & \log(1/2) + 2Y_n n^{-\alpha} -2Y_n^2 n^{-2\alpha}
+ O(n^{-3\alpha}) ,\eean
for all $n$ sufficiently large, and
\bean 
\log q_n & = & 
\log (1/2) + \log ( 1 - 2Y_n n^{-\alpha} ) \\
& = & \log(1/2) - 2Y_n n^{-\alpha} -2Y_n^2 n^{-2\alpha}
+ O(n^{-3\alpha}) ,\eean
so that
\bea
\label{0601s}
\log (p_n/q_n)  = \log p_n - \log q_n =
 4Y_n n^{-\alpha} + O(n^{-3\alpha}).\eea

Lemma \ref{lem0920} below gives
 bounds for the expected hitting time $T(n)$, and so prepares us for the proof of
Theorem \ref{thm0} via an application of Lemma \ref{1002e}.
\begin{lemma}
\label{lem0920} Suppose $\bbP[\xi_1=1/2]=1$, $Y_1 \eqd -Y_1$, 
$\sigma^2 \in (0, \infty)$, and $\alpha \in (0,1/2)$. Then
for a.e.~$\omega$, for any $\eps>0$, 
for all but finitely many $n$,
\bea
\label{0920e}
 \exp ( n^{(1-2\alpha)/2} (\log n)^{-1} (\log \log n)^{-1-\eps})
 \leq
T(n) \nonumber \\
\leq
 \exp ( n^{(1-2\alpha)/2} (\log \log n)^{(1/2)+\eps} ).\eea
\end{lemma}
\proof 
From (\ref{0608a}), $\bbE[\log ( p_k / q_k )]=0$ and from
(\ref{0601s}) $\Var[\log ( p_k / q_k )]=16\sigma^2 k^{-2\alpha} +
o(k^{-2\alpha})$. Hence, for $\alpha \in (0,1/2)$, for all $i$,
\bea
\label{1002f}
 C_1 i^{1-2\alpha} \leq
  \Var \sum_{k=1}^i \log ( p_k / q_k )
\leq C_2 i^{1-2\alpha},\eea
for some $C_1, C_2 \in (0,\infty)$ with $C_1 < C_2$.

Now we derive the lower bound in (\ref{0920e}).
By Lemma \ref{lem0601} and (\ref{1002f}), 
 for an appropriate choice of $a(\cdot)$ satisfying the
 conditions of Lemma \ref{lem0601}, for a.e.~$\omega$, a.s.,
  \bea
  \label{ab3}
   \max_{1 \leq i \leq n-1} \sum_{k=1}^i \log(p_k/q_k) \geq
 C n^{(1-2\alpha)/2} a( n^{1-2\alpha}),\eea
 for all but finitely many $n$.
For $\eps>0$, we take
$a(n) = ( \log n )^{-1}(\log \log n)^{-1-\eps}$; 
then $a(\cdot)$ satisfies the conditions of Lemma 4.
Then (\ref{ffff}) and (\ref{ab3})
imply the
lower bound in (\ref{0920e}).

Now we prove the upper bound in
(\ref{0920e}), using (\ref{eeee}). 
By Lemma \ref{itlog} with (\ref{1002f}) we have that for
a.e.~$\omega$, a.s.,
for all but finitely many $n$,
\bean
 \max_{0 \leq i \leq n-1} \sum_{k=1}^i \log (p_k/q_k)   < C n^{(1-2\alpha)/2}
  (\log \log{n})^{1/2}, \\
  \max_{0 \leq i \leq n-1} \sum_{k=1}^i (-\log (p_k/q_k))   < C n^{(1-2\alpha)/2}
  (\log \log{n})^{1/2},
\eean
 for some $C\in(0,\infty)$. Thus from (\ref{eeee})
 we obtain the upper bound in (\ref{0920e}).
  $\square$\\

\noindent {\bf Proof of Theorem.} First we prove part
(i) of Theorem \ref{thm0}. From the lower bound in
(\ref{0920e}), we have that, for a.e.~$\omega$,
there exists a finite positive constant $C$ (depending on $\omega$)
such that, for any
$\eps>0$, for all $n$ sufficiently large,
\bea
 \label{0601h} T(n) \geq g(n) := C \exp \left( n^{(1-2\alpha)/2}
( \log n )^{-1}(\log \log n)^{-1-\eps}
  \right). \eea
So by (\ref{001}), we have that, for a.e.~$\omega$, for
  any $\eps>0$, a.s.,
  \[ \eta_t (\omega) \leq  g^{-1} ( 4t^2) \leq C
  ((\log t) (\log \log t)^{1+\eps})^{2/(1-2\alpha)},  \]
  for all but finitely many $t$, which gives
(\ref{0210a}).  Now we prove part (ii). From the upper bound in
(\ref{0920e}), we have that, for any $\eps>0$,
\[ T(n) \leq h(n) := C \exp
( n^{(1-2\alpha)/2} (\log \log
n)^{(1/2)+\eps}),
\]
so that, for all $n$ sufficiently large,
\bea
\label{uuu}
 h^{-1} (n) \geq C (\log n)^{2/(1-2\alpha)} (\log \log \log n)^{-(1+3\eps)/(1-2\alpha)}.\eea
From  (\ref{003}) we have that a.s., for infinitely many $t$,
\[ h(\eta_t(\omega)) \geq t (\eta_t(\omega))^{-2}
\geq C t (\log t)^{-5/(1-2\alpha)},\]
by (\ref{0210a}).  Thus a.s., for infinitely many $t$,
$\eta_t(\omega) \geq h^{-1} (
C t (\log t)^{-5/(1-2\alpha)} )$,
which with  (\ref{uuu}) yields
(\ref{0920c}). 
 $\square$

\subsection{Proofs of Theorems \ref{thm5} and \ref{thm8}}
\label{secprf2}

To prove Theorems \ref{thm5} and \ref{thm8}, we proceed along the
same lines as the proof of Theorem \ref{thm0} in Section
\ref{secprf1}, and apply Lemma \ref{1002e}. Theorem \ref{thm8} (the
transient case) requires some extra work, both to obtain
suitable bounds for $T(n)$ and to prove
that the lower bound on the random walk
holds
 all but finitely often.

Suppose $\bbE[\zeta_1]=0$, $s^2 \in (0,\infty)$,
 $\sigma^2 \in [0,\infty)$.  Then for
 a.e.~$\omega$
\bea
\label{0708q}
 \log \left( \frac{p_n}{q_n} \right)
= \zeta_n + \log \left( 1+ \frac{Y_n}{\xi_n} n^{-\alpha}\right) -
\log \left( 1 - \frac{Y_n}{1-\xi_n} n^{-\alpha}\right),\eea  for all $n \geq n_0$ for a finite 
 absolute constant $n_0$, where
$\zeta_i$, $i\in\N$, as defined at (\ref{0520b}) are i.i.d.~with $\bbE[\zeta_1]=0$
and $\Var[\zeta_1] \in (0,\infty)$.
 It follows from (\ref{0708q}) and Taylor's theorem
  that, for a.e.~$\omega$, for all $n$
sufficiently large, 
\bea \label{0427f} \log ( p_n / q_n
) = \zeta_n + Z_n n^{-\alpha}+O(n^{-2\alpha}),\eea where
$Z_i$, $i\in\N$, are i.i.d.~with $\bbE[Z_1]=\lambda$ (see
(\ref{0520b}) and (\ref{0520bx})).
 Then by (\ref{0427f})
\bea
\label{0439a}
 \bbE[\log(p_k/q_k)] = \lambda k^{-\alpha} + O(k^{-2\alpha}), 
 \Var[\log(p_k/q_k)] = \Var[ \zeta_1 ] + O(k^{-\alpha}) .\eea

\begin{lemma}
\label{lem0430a} Suppose $\bbE[\zeta_1]=0$, $s^2 \in (0,\infty)$,
 $\sigma^2 \in [0, \infty)$, $\lambda<0$, and $\alpha \in (0,1/2)$. 
 For a.e.~$\omega$ and any $\eps>0$, 
 for all but finitely many $n$,
  \bea \label{0427e}
 \exp ( n^\alpha (\log n)^{-2-\eps} ) \leq
T(n) \leq
 \exp ( n^\alpha ( \log n)^{1+\eps} ).\eea
 \end{lemma}
\proof First we prove the upper bound in (\ref{0427e}). 
Since
$\lambda<0$, we have from (\ref{0439a}) that
\bea
\label{tt1}
 \bbE \sum_{k=i-j+1}^i \log(p_k/q_k) 
 \leq - C ( i^{1-\alpha} - (i-j)^{1-\alpha} )
 ,\eea
for some $C \in  (0,\infty)$. Taylor's theorem implies that for $\alpha \in (0,1)$
\bea
\label{tt2}
 i^{1-\alpha} - (i-j)^{1-\alpha} = C j i^{-\alpha} (1- \theta (j/i) )^{-\alpha} ,\eea
for some $C \in (0,\infty)$ and $\theta \in (0,1)$. Thus 
it follows from (\ref{tt1}) and (\ref{tt2}) 
 that 
for all $i \in \N$, and all $j =1,2,\ldots, i$
\bea
\label{0429q}
 \bbE \sum_{k=i-j+1}^i \log(p_k/q_k) \leq -C j i^{-\alpha} ,
\eea
for some $C\in (0,\infty)$. 
By Lemma \ref{ubd1}  we
have that, for some $C \in (0,\infty)$, for a.e.~$\omega$, all but finitely many $i$,
and all $j=1,2,\ldots,i$,
\bea
\label{0429p}
   \sum_{k=i-j+1}^i (\log(p_k/q_k)-\bbE[\log(p_k/q_k)])  \leq C j^{1/2} (\log i)^{1/2}.
\eea 
 Suppose $\eps>0$. Then from (\ref{0429p}) with (\ref{0429q}), for a.e.~$\omega$,
for $j \geq \lceil i^{2\alpha} (\log i)^{1 + \eps}\rceil$ \bea
\label{0430a}
 \sum_{k=i-j+1}^i \log(p_k/q_k) \leq -C ji^{-\alpha}+C' j^{1/2} (\log i)^{1/2}
 \leq -C'' ji^{-\alpha},\eea
and, for $j \leq \lceil i^{2\alpha} (\log i)^{1 + \eps}\rceil$ \bea
\label{0430b}
 \sum_{k=i-j+1}^i \log(p_k/q_k) \leq C j^{1/2} (\log i )^{1/2},\eea
where each inequality holds for all but
finitely many $i$. So from (\ref{1001c}), (\ref{0430a}) and (\ref{0430b}) we obtain, for
a.e.~$\omega$, for any $\eps>0$, for all but finitely many $i$,
\bean \Delta_i & \leq & \sum_{j=0}^{\lceil i^{2\alpha} (\log
i)^{1+\eps} \rceil} \exp (C j^{1/2} ( \log i)^{1/2})
+ \sum_{j=\lceil i^{2\alpha} ( \log{i})^{1+\eps}\rceil}^i \exp( -C' j i^{-\alpha}) \\
& \leq & \exp ( C'' i^{\alpha} (\log i)^{1+\eps} ),\eean 
for $C'' \in (0,\infty)$. Then the
upper bound for $T(n)$ in (\ref{0427e}) follows.

We now prove the lower bound in (\ref{0427e}). For $\eps>0$ set $k_\eps (1):=1$ and for
 $i
>1$ define
\bea \label{1002w} k_\eps (i) := \lfloor i^{2 \alpha} (\log
i)^{-2-\eps} \rfloor .\eea 
Then, for any $\alpha \in (0,1/2]$ and all $n$ sufficiently large,
from (\ref{ab1}),
\bea
\label{45a}
 T(n) \geq \max_{1\leq i \leq n-1} \max_{1 \leq j \leq k_\eps (i) }
\exp \sum_{k=i-j+1}^i \log (p_k/q_k).\eea
Then (\ref{0708q}) and
Taylor's theorem imply that there is a constant $C\in (0,\infty)$ such that,
for all $k$, $\log(p_k/q_k) = \zeta_k + W_k k^{-\alpha}$,
where $|W_k| < C$. Thus for $i \in \N$ and $j=1,2,\ldots,i$,
\[ \sum_{k=i-j+1}^i \log (p_k/q_k)
\geq \sum_{k=i-j+1}^i \zeta_k - C \sum_{k=i-j+1}^i k^{-\alpha}
\geq \sum_{k=i-j+1}^i \zeta_k - C' j i^{-\alpha},\]
again using Taylor's theorem (cf (\ref{tt2})). Hence by (\ref{45a})
\bea
\label{pp1}
 T(n)\geq \exp \left(
\max_{1\leq i \leq n-1} \max_{1 \leq j \leq k_\eps (i) }
\sum_{k=i-j+1}^i \zeta_k - C k_\eps (n) n^{-\alpha}\right).\eea
By Lemma \ref{lbd1}, we have that for any $\eps>0$,
for a.e.~$\omega$,
\[ \max_{1\leq i \leq n-1} \max_{1 \leq j \leq k_\eps (i) }
\sum_{k=i-j+1}^i \zeta_k \geq ( k_\eps (n/2))^{1/2} ( \log n)^{-1-(\eps/4)}
\geq C n^{\alpha} (\log n)^{-2-(3\eps/4)},\]
for all but finitely many $n$, while
$k_\eps (n) n^{-\alpha} \leq n^{\alpha} (\log n)^{-2-\eps}$.
Hence (\ref{pp1}) implies the lower bound in (\ref{0427e}). $\square$

\begin{lemma}
 Suppose $\bbE[\zeta_1]=0$, $s^2 \in (0,\infty)$,
 $\sigma^2 \in[0,\infty)$, and $\lambda \neq
0$. 
\begin{itemize}
\item[(i)] Suppose  $\alpha > 1/2$. For a.e.~$\omega$ and any $\eps>0$, 
for all but finitely many $n$,
\bea \label{0920d} \exp ( n^{1/2} (\log n)^{-1-\eps})
\leq T(n) \leq
 \exp ( n^{1/2} (\log \log n)^{(1/2)+\eps} ).\eea
\item[(ii)] Suppose  
$\alpha = 1/2$. For a.e.~$\omega$ and any $\eps>0$, 
for all but finitely many $n$,
\bea \label{0920dd} \exp ( n^{1/2} (\log n)^{-2-\eps})
\leq T(n) \leq
 \exp ( n^{1/2} (\log \log n)^{(1/2)+\eps} ).\eea
\end{itemize}
\end{lemma}
\proof
To prove the upper bounds in (\ref{0920d}) and (\ref{0920dd}), we apply
(\ref{eeee}).
For $\lambda \neq 0$, $\alpha \geq 1/2$ 
 we have from (\ref{0439a}) that 
 \[\sum_{k=1}^i \bbE[\log(p_k/q_k)] = O( \max \{ i^{1-\alpha} , \log i \} ),\]
 so that for some $C \in (0,\infty)$ and all $n$
\bea
\label{ccc1}
 \max_{0 \leq i \leq n} \sum_{k=1}^i \log (p_k/q_k)
\leq \max_{0 \leq i \leq n} \sum_{k=1}^i (\log (p_k/q_k)
-\bbE[\log(p_k/q_k)]) +C \max \{ n^{1-\alpha} , \log n \} ;\eea
similarly for the second maximum in (\ref{eeee}). By Lemma \ref{itlog}
and (\ref{0439a}),
for a.e.~$\omega$,
\[  \max_{0 \leq i \leq n} \sum_{k=1}^i (\log (p_k/q_k)
-\bbE[\log(p_k/q_k)]) \leq C n^{1/2} (\log\log n)^{1/2} ,\]
for all but finitely many $n$, and since $\alpha \geq 1/2$, (\ref{ccc1})
then implies that for a.e.~$\omega$,
\[ \max_{0 \leq i\leq n} \sum_{k=1}^i \log (p_k/q_k) \leq C n^{1/2} (\log \log n)^{1/2},\]
 for all but finitely many $n$, 
and similarly for the second maximum in (\ref{eeee}).
Then (\ref{eeee}) gives the upper bounds in (\ref{0920d}) and (\ref{0920dd}).

Now we prove the lower bounds in (\ref{0920d}) and (\ref{0920dd}).
 In the case $\alpha > 1/2$,
 \[ \max_{1 \leq i \leq n-1} \sum_{k=1}^i \log (p_k/q_k)
 \geq  \max_{1 \leq i \leq n-1} \sum_{k=1}^i (\log (p_k/q_k)
 -\bbE[ \log (p_k/q_k) ])
 - C\max \{ n^{1-\alpha} , \log n \} ,\]
 by a similar argument to (\ref{ccc1}). 
 Lemma \ref{lem0601} implies that for any $\eps>0$, for a.e.~$\omega$,
 \[ \max_{1 \leq i \leq n-1} \sum_{k=1}^i (\log (p_k/q_k)
 -\bbE[ \log (p_k/q_k) ]) \geq n^{1/2} (\log n)^{-1-\eps},\]
 for all but finitely many $n$;
then (\ref{ffff}) implies the lower bound in
 (\ref{0920d}). Finally, suppose $\alpha=1/2$. Once more we define
 $k_\eps(i)$ by (\ref{1002w}), and follow the argument for
 (\ref{pp1}). This yields the lower bound in (\ref{0920dd}).
$\square$\\

\noindent {\bf Proof of Theorem \ref{thm8}.} 
For the upper bound in (\ref{0427c}), the lower bound  in
(\ref{0427e}) implies that, for a.e.~$\omega$, for any $\eps>0$
there exists $C\in (0,\infty)$ such that
\[ T(n) \geq g(n) := C \exp (n^\alpha (\log n)^{-2-\eps}),\]
for all $n$ sufficiently large.
 Then (\ref{001}) gives, for a.e.~$\omega$, for
any $\eps>0$, a.s.,
\[ \eta_t(\omega) \leq  g^{-1} (4t^2) \leq C (\log t)^{1/\alpha} (\log
\log t)^{(2+2\eps)/\alpha},  \] for all but finitely many $t$. Then
the upper bound in (\ref{0427c}) follows.

We now want to obtain the lower
bound in (\ref{0427c}). Recalling the proof of Lemma \ref{1002e}(ii),
we were able to show that, along a sequence of first hitting times for the
random walk,
these  times were not too large. This gave us a lower bound that was valid infinitely often.
In order to extend this technique to the transient case, and obtain a lower bound
valid {\em all but finitely} often, we show in addition that
(roughly speaking),
 in the present case, the time of the last visit of the random walk
to a site is not too much greater than the first hitting time.
 
 For fixed $\omega$, let $a_n$ denote the probability that the random walk $\eta_t(\omega)$
 hits $n$ in finite time, given that it starts at $2n$. 
 For $n \geq 1$ define
 \bea
 \label{ooo1}
  M_n := 1+ \sum_{j=1}^\infty \prod_{k=1}^{j} \frac{p_{n+k}}{q_{n+k}}
 = 1+ \sum_{j=1}^\infty \exp \sum_{k=1}^{j} \log \left( \frac{p_{n+k}}{q_{n+k}}\right)
 .\eea
 Standard hitting 
 probability arguments yield $a_0=1$, and for $n \geq 1$,
 if $M_n < \infty$,
 \bea
 \label{ooo2}
  a_n = M_n^{-1} \sum_{j=n}^\infty \prod_{k=1}^{j} \frac{p_{n+k}}{q_{n+k}}
 = M_n^{-1} \sum_{j=n}^\infty \exp \sum_{k=1}^{j} \log \left( \frac{p_{n+k}}{q_{n+k}}\right).
 \eea
In the present case ($\lambda<0$, $\alpha \in (0,1/2)$), (\ref{0439a}) holds. Thus for $n, j \in \N$
\bea
\label{ss1}
\bbE \sum_{k=1}^j \log (p_{n+k}/q_{n+k}) \leq -C ( (n+j)^{1-\alpha} - n^{1-\alpha} ) 
,\eea
for some $C \in (0,\infty)$.
Here, by Taylor's theorem, for $\alpha \in (0,1)$,
\bea
\label{ss3}
 (n+j)^{1-\alpha} - n^{1-\alpha} = C j (n+\theta j)^{-\alpha} ,\eea
for some $C \in (0,\infty)$, $\theta \in (0,1)$. In particular, for
$j \geq n$, (\ref{ss1}) and (\ref{ss3}) imply
\bea
\label{ss2} \bbE \sum_{k=1}^j \log (p_{n+k}/q_{n+k})\leq - Cj (j (1+\theta))^{-\alpha}
\leq -C' j^{1-\alpha},\eea
for $C' \in (0,\infty)$. 
Also, by the Azuma-Hoeffding inequality and an argument similar to Lemma 
\ref{ubd1}, we have that, for a.e.~$\omega$,
\[ \sum_{k=1}^j (\log(p_{n+k}/q_{n+k}) -\bbE [ \log(p_{n+k}/q_{n+k})])
\leq C j^{1/2} (\log (jn) )^{1/2} ,\]
for all but finitely many $(n,j)$. Thus for all $(n,j)$ we have  that, for a.e.~$\omega$,
\bea
\label{ss5}
 \sum_{k=1}^j \log(p_{n+k}/q_{n+k}) \leq C j^{1/2} (\log (jn) )^{1/2},\eea
 for some $C \in (0,\infty)$.
However,  we have
from  (\ref{ss2}) that, for a.e.~$\omega$, there are constants $C,C',C'' \in (0,\infty)$
such that, for all   $n\in \N$, and $j \geq n$
\bea
\label{ss4}
 \sum_{k=1}^j \log(p_{n+k}/q_{n+k}) \leq - C j^{1-\alpha} + C' j^{1/2} ( \log {j} )^{1/2}
\leq -C'' j^{1-\alpha},\eea
  since $\alpha \in (0,1/2)$.
Hence, for a.e.~$\omega$,  
from (\ref{ooo1}), (\ref{ss5}), and (\ref{ss4}), for $n \in \N$,
\bean M_n \leq   \sum_{j=1}^{n}
\exp ( C j^{1/2} (\log (jn) )^{1/2} )
+ \sum_{j=n}^\infty
\exp ( -C' j^{1-\alpha} ) \\
\leq \exp ( C'' n^{1/2} (\log n)^{1/2} ) < \infty.\eean
Further, since $M_n \geq 1$ for all $n$,
  (\ref{ooo2}) and (\ref{ss4}) imply,
for a.e.~$\omega$, for  all $n\in \N$,
\[ a_n \leq  \sum_{j=n}^\infty \exp ( -C j^{1-\alpha}  )
\leq \exp ( -C' n^{1-\alpha}  ),\]
 for some $C' \in (0,\infty)$. Thus, for a.e.~$\omega$, 
 $\sum_n a_n <\infty$.

The (first)
Borel-Cantelli
lemma then implies that, for a.e.~$\omega$, a.s.,
 for only finitely many sites $n$ 
does $\eta_t (\omega)$ return to $n$ after visiting $2n$. Denoting by $\ell_n$ the 
time of the
last
visit of $\eta_t(\omega)$ to $n$, we then have that $\ell_n \leq \tau_{0,2n}$  a.s.~for all
but finitely many $n$. Suppose $T(n) \leq h(n)$ for all $n$.
Following the proof of Lemma \ref{1002e}(ii), we 
have that  a.s., for all but finitely many $n$, $\tau_{0,n} \leq n^2 h(n)$. Thus  
 for a.e.~$\omega$, a.s.,
 \bea
 \label{bbx}
  \ell_n \leq \tau_{0,2n} \leq 4n^2 h(2n),\eea
for all
$n \geq n_0$ for some  finite $n_0$ (depending on $\omega$).

Moreover, since, for a.e.~$\omega$, $\eta_t(\omega)$
is transient, we have that, for a.e.~$\omega$,  a.s., $\eta_t(\omega) \geq n_0$ for all $t$ sufficiently
large. Hence from (\ref{bbx}), using the fact that $\ell_{\eta_t(\omega)} \geq t$ for all $t$, we have
that for a.e.~$\omega$, a.s., for all but finitely many $t$, 
\[ 
t \leq \ell_{\eta_t(\omega)}  \leq  4 \eta_t(\omega)^2 h(2 \eta_t(\omega)).  \]
Then,  with the upper bound
  in (\ref{0427e}), we obtain, for a.e.~$\omega$, for any
$\eps>0$,  a.s., 
\[ t < \exp  ( \eta_t(\omega)^\alpha
(\log  \eta_t(\omega))^{1+\eps}  ) ,\]
for all but finitely many $t$. This implies the
lower bound in (\ref{0427c}).
$\square$\\

\noindent {\bf Proof of Theorem \ref{thm5}.} We first prove part
(i). Suppose $\alpha>1/2$. From the lower bound on
$T(n)$ in (\ref{0920d}), for a.e.~$\omega$, for any $\eps>0$,
 \bean T(n)
\geq g(n) := C \exp (n^{1/2} (\log n)^{-1-\eps}),\eean for all  
$n \in \N$. Then (\ref{001}) implies the upper bound in (\ref{0210ac}).
For part (ii), when $\alpha=1/2$, the lower bound   in
(\ref{0920dd}) allows us, this time, to take $g(n) := C \exp (n^{1/2} (\log n)^{-2-\eps})$.
Then (\ref{001}) gives the upper bound in (\ref{0210acd}).

For part (iii) of the theorem, for $\alpha \geq
1/2$, the upper bound on $T(n)$ in (\ref{0920d}) and
(\ref{0920dd}) implies that for a.e.~$\omega$ 
\[T(n) \leq h(n) := C\exp ( n^{1/2} (\log\log n)^{(1/2)+\eps} ),\]
for all but finitely many $n$; in particular $h^{-1}$ satisfies the lower bound
of  (\ref{uuu}) with $\alpha =0$.  Then (\ref{003}) yields
the lower bound in (\ref{0920a}). 
$\square$

\subsection{Proofs of Theorems \ref{thm21} and \ref{thm20}}
\label{secprf2a}
 
We now move on to the ergodic cases (Theorems \ref{thm21} and \ref{thm20}).
Again  we start
by bounding $T(n)$. First we deal with the ergodic case of the 
random perturbation of the simple random walk.

\begin{lemma} 
\label{lem33}
Suppose $\bbP[\xi_1=1/2]=1$, $\bbE[Y_1]>0$, $\sigma^2 \in(0,\infty)$, and
$\alpha \in (0,1)$. Then for a.e.~$\omega$, as $n \to \infty$
\bean
T(n) = \exp \left( \frac{4 \bbE[Y_1]}{1-\alpha} n^{1-\alpha} [1+o(1)] \right).
\eean
\end{lemma}
\proof  In this case, (\ref{0601s}) holds. We apply a variation of the argument 
for Lemma \ref{ubd1}. We have that for each $i$
\[ Y^i_j := \sum_{k=i-j+1}^i ( \log(p_k/q_k)- \bbE[\log(p_k/q_k)])\]
is a martingale over $j=1,2,\ldots,i$, with increments $|Y^i_{j}-Y^i_{j-1} |$
bounded by
\[  |\log(p_{i-j+1}/q_{i-j+1}) |
+ | \bbE[\log(p_{i-j+1}/q_{i-j+1})] | \leq C (i-j+1)^{-\alpha} =: c_j^i ,\]
for some $C \in (0,\infty)$,
by (\ref{0601s}). Thus for each $j \leq i$, for $\alpha \in (0,1)$,
\[ \sum_{k=1}^j (c_k^i)^2 =  C \sum_{k=1}^j (i-k+1)^{-2\alpha} \leq C' i^{1-\alpha}.\]
 Then for each $i$ and $j \leq i$
  the Azuma-Hoeffding inequality implies that 
\[ \bbP [ |Y^i_j |\geq t ]\leq 2 \exp ( - C t^2 i^{\alpha-1} ),\]
for all $t>0$.
Hence for any $\eps>0$, the Borel-Cantelli lemma implies that
\[ \max_{1 \leq j \leq i}
|Y_j^i | \leq i^{((1-\alpha)/2)+\eps},\]
for all but finitely many $i$. Also, from (\ref{0601s}),
\[ \bbE \sum_{k=i-j+1}^i \log (p_k/q_k) = \frac{4 \bbE[Y_1]}{1-\alpha} \left( i^{1-\alpha}
- (i-j)^{1-\alpha} \right) [1+o(1)].\] Hence for all $i$ sufficiently large,
since $\eps>0$ was arbitrary and $\alpha\in (0,1)$
\bea
\label{gg1}
 \sum_{k=i-j+1}^i \log (p_k/q_k) = \frac{4 \bbE[Y_1]}{1-\alpha} \left( i^{1-\alpha}
- (i-j)^{1-\alpha} \right)[1+o(1)] + o(i^{1-\alpha}).\eea
Thus from (\ref{1001c}) and (\ref{gg1}),
 as $i\to\infty$,
\bean
\Delta_i & = &   \exp \left( \frac{4 \bbE[Y_1]}{1-\alpha} i^{1-\alpha}
 [1+o(1)]  \right) \sum_{j=0}^i \exp \left( -Cj^{1-\alpha}[1+o(1)] \right)
  \\
 & = & \exp \left( \frac{4 \bbE[Y_1]}{1-\alpha} i^{1-\alpha}
 [1+o(1)]  \right),\eean
 from which the lemma follows. 
$\square$\\

\noindent {\bf Proof of Theorem \ref{thm20}.} Once again we apply Lemma \ref{1002e}.
First we prove
the lower bound. From Lemma \ref{lem33} we have that
for a.e.~$\omega$, for all $n$ 
\[ T(n) \geq g(n) := \exp \left( \frac{4 \bbE[Y_1]}{1-\alpha} n^{1-\alpha} [1+o(1)] \right).\]
It follows that
\[ g^{-1} (n) = \left( \frac{1-\alpha}{4 \bbE[Y_1]} \right)^{1/(1-\alpha)}
(\log n)^{1/(1-\alpha)} [1+o(1)].\]
Then (\ref{001}) implies that a.s., for all but finitely many $t$, 
for any $\eps>0$
\bean \eta_t(\omega) \leq g^{-1} ( (2t)^{1+\eps}) \\
=
(1+\eps)^{1/(1-\alpha)} \left( \frac{1-\alpha}{4 \bbE[Y_1]} \right)^{1/(1-\alpha)}
(\log t)^{1/(1-\alpha)} [1+o(1)],\eean
and thus we obtain the upper bound in the theorem.
 On the other hand,
Lemma \ref{lem33} implies that
for a.e.~$\omega$, any $\eps>0$, and all  $n$  
\[ T(n) \leq  h(n) := \exp \left( \frac{4 \bbE[Y_1]}{1-\alpha} n^{1-\alpha}  [1+o(1)]  \right).\]
Then  (\ref{003}) implies that, a.s., for infinitely many $t$,
\[ t \leq (\eta_t (\omega))^2 h( \eta_t(\omega))
\leq \exp \left( \frac{4 \bbE[Y_1]}{1-\alpha} (\eta_t (\omega))^{1-\alpha}  [1+o(1)]  \right),\]
  from which the lower bound in the theorem follows. 
$\square$\\

Now we deal with the ergodic case of the perturbation of Sinai's regime.

\begin{lemma} 
\label{lem34}
Suppose $\bbE[\zeta_1]=0$, $ s^2 \in (0,\infty)$, $\lambda>0$, $\sigma^2 \in [0,\infty)$, and
$\alpha \in (0,1/2)$. Then for a.e.~$\omega$, as $n \to \infty$
\bean
T(n) = \exp \left( \frac{\lambda}{1-\alpha} n^{1-\alpha} [1+o(1)] \right).
\eean
\end{lemma}
\proof In this case, we have that (\ref{0439a}) holds (now with $\lambda>0$). Thus
\[ \bbE \sum_{k=i-j+1}^i \log (p_k /q_k) = \frac{\lambda}{1-\alpha} \left( i^{1-\alpha}
-(i-j)^{1-\alpha} \right) [1+o(1)] .\]
Now we can apply Lemma \ref{ubd1} 
to obtain for a.e.~$\omega$, for all but finitely many $i$,
\[ \left| \sum_{k=i-j+1}^i (\log (p_k /q_k)
-\bbE[ \log (p_k /q_k)]) \right| \leq C j^{1/2} (\log i)^{1/2},\]
for $j=1,2,\ldots,i$.
Since $\alpha<1/2$ 
we have that for a.e.~$\omega$, as $i \to \infty$
\bea
\label{gg2} 
\sum_{k=i-j+1}^i \log (p_k /q_k) = \frac{\lambda}{1-\alpha} \left( i^{1-\alpha}
-(i-j)^{1-\alpha} \right) [1+o(1)] .\eea
Hence, from
(\ref{1001c}) and (\ref{gg2}), as $i \to \infty$,
\bean
\Delta_i & = & \exp \left(\frac{\lambda}{1-\alpha} i^{1-\alpha}
[1+o(1)] \right) \sum_{j=0}^i \exp \left( -C j^{1-\alpha} [1+o(1)]  \right)\\
& = & \exp \left(\frac{\lambda}{1-\alpha} i^{1-\alpha}
[1+o(1)] \right),
\eean
and so the lemma follows.  $\square$\\

\noindent {\bf Proof of Theorem \ref{thm21}.} The proof follows in a similar
way to the above proof of Theorem \ref{thm20}, this time using the bounds in Lemma \ref{lem34}
and applying Lemma \ref{1002e} once more. $\square$

\subsection{Proof of Theorem \ref{thm3}}
\label{prf3}

We now prove Theorem \ref{thm3}.
Once more, with the definition of
$\zeta_n$ and $Z_n$ at (\ref{0520b}), we have that for a.e.~$\omega$ and $n$ sufficiently large,
$\log(p_n/q_n)$
is given by (\ref{0708q}), (\ref{0427f}).
 In this
case $\bbE[\zeta_1]=0$ and $Y_1/\xi_1 \eqd -Y_1/(1-\xi_1)$, which 
 implies that $\bbE[\log(p_n/q_n)]=0$ for all
$n$ sufficiently large.

\begin{lemma}
\label{lem0920q} Suppose $\bbE[\zeta_1]=0$, $s^2 \in (0,\infty)$,
$Y_1/\xi_1 \eqd -Y_1/(1-\xi_1)$, $\sigma^2\in[0, \infty)$, and
$\alpha>0$. For a.e.~$\omega$, for any $\eps>0$,
for all but finitely many $n$,
\bea \label{0920g} \exp( n^{1/2} (\log n)^{-1-\eps} ) \leq
 T(n)
\leq
 \exp ( n^{1/2} (\log \log n)^{(1/2)+\eps} ).\eea
\end{lemma}
\proof We apply Lemma \ref{lowerbd}. 
We have
that (\ref{0427f}) and (\ref{0439a}) hold in this case. 
For the upper bound, consider (\ref{eeee}).
By Lemma
\ref{itlog} we have that for a.e.~$\omega$, for all but
finitely many $n$, \bean
 \max_{0 \leq i \leq n-1}
 \sum_{k=1}^i \log (p_k/q_k)  < C n^{1/2} (\log \log n )^{1/2},
\eean for some $C\in(0,\infty)$, and similarly
for the second maximum in (\ref{eeee}). 
Then (\ref{eeee}) implies
the upper bound in (\ref{0920g}). For the lower
 bound, we use (\ref{ffff}).
We apply Lemma \ref{lem0601} with $a(x)= (\log x)^{-1-\eps}$ to obtain, for a.e.~$\omega$,
 for any $\eps>0$
\[ \max_{1 \leq i \leq n-1}
 \sum_{k=1}^i \log (p_j/q_j) \geq n^{1/2} (\log n)^{-1-\eps},\]
 for all but finitely many $n$. With (\ref{ffff}), the lower bound in (\ref{0920g})
 follows.
 $\square$\\

\noindent {\bf Proof of Theorem.} Again the proof is very similar
to that of Theorems \ref{thm0} and \ref{thm5}, this time using
 Lemma \ref{lem0920q} and Lemma \ref{1002e}. $\square$

\subsection{Proofs of Theorems \ref{thm9} and \ref{thm10}}

Finally, we prove the results on the stationary distribution in the ergodic cases given in
Theorems \ref{thm9} and \ref{thm10}. Given $\omega$, suppose $\eta_t(\omega)$ is ergodic; then there exists a unique stationary distribution $(\pi_0, \pi_1, \pi_2, \ldots)$. It is straightforward to obtain the result (see, for example, Lemma 5 of \cite{mw}) that,
for a given    $\omega$ such that $\eta_t(\omega)$ is ergodic,
there exists a constant $C\in(0,\infty)$ such that, for all $n \geq 2$,
\bea
\label{0922e}
\pi_n = C \prod_{k=1}^{n} \frac{q_k}{p_k}
= C \exp \left( \sum_{k=1}^n \log (q_k/p_k) \right).
\eea

\noindent
{\bf Proof of Theorem \ref{thm9}.} Here we have that $\log(p_n/q_n)$ is given by
(\ref{0427f}), with  $\lambda>0$, $\alpha \in (0,1/2)$ and $\bbE[\zeta_1]=0$.
In this case the $j=i=n$ case of (\ref{gg2}) implies that
\[ \sum_{k=1}^n \log (q_k/p_k) = -\sum_{k=1}^n \log (p_k/q_k)
= - \frac{\lambda}{1-\alpha} n^{1-\alpha} [1+o(1)],\]
as $n\to \infty$.
Then (\ref{0922e}) yields (\ref{0922a}). $\square$\\

\noindent
{\bf Proof of Theorem \ref{thm10}.} This time we have that $\log(p_n/q_n)$ is given by
(\ref{0601s}), where now $\bbE[Y_1]>0$ and $\alpha \in (0,1)$. 
In this case the $j=i=n$ case of (\ref{gg1}) implies that
\[ \sum_{k=1}^n \log (q_k/p_k) =- \sum_{k=1}^n \log (p_k/q_k)
= - \frac{4\bbE[Y_1]}{1-\alpha} n^{1-\alpha} [1+o(1)],\]
as $n \to \infty$. Then (\ref{0922e}) yields (\ref{0922c}). $\square$

\begin{center} \textbf{Acknowledgements} \end{center}
  Some of this work was done when AW was at the University of Durham, supported by
an EPSRC doctoral training account, and subsequently at the University of Bath.
We are grateful to Serguei Popov for useful discussions,
and to
an anonymous referee for a careful reading of an earlier version of this paper.

\end{document}